\documentstyle{amsppt}
\overfullrule=0pt

\magnification=\magstep1
 \pagewidth{6.0 true in}
\pageheight{9.0 true in}

\hoffset 0in % change for fine tuning
\voffset -0.10in

\parindent 0pt
\document

\topmatter
\title Remarks on weakly pseudoconvex boundaries: Erratum
\endtitle
\author Judith Brinkschulte \\
C. Denson Hill \\
Mauro Nacinovich
\endauthor
\address Judith Brinkschulte - Universit\"at Leipzig, Mathematisches Institut,
Augustusplatz 10/11, 04109 Leipzig. Germany\endaddress
\email brinkschulte\@math.uni-leipzig.de \endemail
\address C.Denson Hill - Department of Mathematics, SUNY at Stony Brook,
Stony Brook NY 11794, USA \endaddress
\email dhill\@math.sunysb.edu \endemail
\address Mauro Nacinovich - Dipartimento di Matematica -
Universit\`a di Roma "Tor Vergata" - via della Ricerca Scientifica -
00133 - Roma - Italy \endaddress
\email nacinovi\@mat.uniroma2.it\endemail
\endtopmatter

\bigskip

It is necessary to make two tiny corrections in  \cite{BHN}. This is because the last two authors, having become aware of some inconsistencies in  \cite{HN1}, have repaired the situation in  \cite{HN2}, but the results they obtained are slightly different than what was originally claimed. Namely, for local results, it is important to make a distinction between the vanishing of the cohomology of small domains and the validity of the Poincar\'e lemma. None of the global results in  \cite{BHN},  concerning weakly pseudoconvex boundaries, require any correction.

\medskip The first change needed is that in Theorem 1, part (ii), on page 2, one must add to the hypotheses that $x_0$ be a {\it regular point} in the sense of  \cite{HN2}.

\medskip 

The second change needed is in the example on page 4. The interesting feature of that example (that it is possible to have {\it vanishing global cohomology}, and at the same time, {\it infinite dimensional local cohomology}) remains true, but it needs to be re-explained. As a bonus we obtain something new and interesting from it.

\medskip
 Let $z=(z_0, z_1)$ be coordinates in $\Bbb C^{2}$, 
$w=(w_1,\ldots ,w_{n-1})$ be coordinates in $\Bbb C^{n-1}$. Consider the egg 
in $\Bbb C^{n+1}$ defined by
$$\Omega = \{ \vert z_0\vert^2 +\vert z_1\vert^2  + \vert w_1\vert^{2m} +\ldots +
\vert w_{n-1}\vert^{2m} < 1\},$$
for an integer $m\geq 2$. It has a weakly pseudoconvex boundary $\partial\Omega$. For $r=0,1,\ldots ,n-1$, 
let $\Sigma_{n-r}$ be the set of points on $\partial\Omega$ at which exactly $r$ 
components of $w$ are zero. Then $\partial\Omega=\bigcup_{k=1}^n \Sigma_k$, and at each 
point $x_0$ of $\Sigma_k$, the Levi form of $\partial\Omega$ has $k$ positive and $n-k$ 
zero eigenvalues. We do not obtain that the Poincar\'e lemma fails at $x_0$ in degree $(p,k)$, as was previously claimed.

\medskip 

However if $B(x_0,r)$ is any sufficiently small ball centered at $x_0\in\Sigma_k$ of radius $r$, in any euclidean metric in $\Bbb C^{n+1}$, it was proved in \cite{HN2}  that
$$\dim H^{p,k}(\partial\Omega\cap B(x_0,r)) = +\infty$$
for all $0\leq p\leq n+1$.
\medskip
 Here is the new observation: Set $U^- =\overline\Omega\cap B(x_0,r)$, $U^+ = \complement\Omega\cap B(x_0,r)$, and $M =\partial\Omega\cap B(x_0,r)$. By \cite{AH1} we have, for the cohomology of smooth forms on the half open/closed domains $U^\pm$, that
$$H^{p,k} (M) = H^{p,k}(U^-) \oplus H^{p,k}(U^+).$$
However it follows from \cite{D} or \cite{N} that $\dim H^{p,k}(U^-) =0$. Thus
$$\dim H^{p,k}(U^+) = +\infty$$
for all $0\leq p\leq n+1$, which is a new result.

\medskip

Since $B(x_0,r)$ is sufficiently small, the Levi form of $M$ has at least $k+1$ positive eigenvalues at each point of $M\setminus\Sigma_k$, but only $k$ positive eigenvalues along $\Sigma_k$. Note that $\Sigma_k$ has real codimension $2n-2k$ in $M$. Now if $\Sigma_k$ had been void, we would know from \cite{AH2}, that it would be possible to choose the Riemannian metric in $\Bbb C^{n+1}$ in such a way as to obtain
$$\dim H^{p,k}(U^+) =0$$
for all $0\leq p\leq n+1$.
Hence we see that the loss of just one positive eigenvalue along the high codimensional locus $\Sigma_k$ in $M$ is enough to convert the cohomology of $U^+$ in degree $k$ from being zero to being infinite dimensional. This was not known before.

\medskip

%%%%%%%%%%%%%%%%%%%%%%%%%%%%%%%%%%%%%%%%%%%%%%%%%%%%%%%%%%%%%%%%%%%%%%%%%%%%%%%%%%%%%%

\end{document}